\documentclass[a4paper,12pt,leqno]{article}

\title{Kropina spaces of constant curvature II.\\ (long version)}
\author{By\\ R. Yoshikawa and K. Okubo}
\date{}
\newlength{\topdummy}
\makeatletter
 
  \@addtoreset{equation}{section}
 \makeatother

\begin{document}
\maketitle
\textbf{Introduction. }\hspace{0.5in}
In 1991, Matsumoto considered the necessary and sufficient conditions for a Kropina space to be of constant curvature 
 and obtained \\

\textbf{Theorem M}([3], [4]) \hspace{0.5in}
\textit{An $n$-dimensional Kropina space is of constant curvature $K$, if and only if, putting $r_{ij}=(b_{i;j}+b_{j;i})/2$
and $s_{ij}=(b_{i;j}-b_{j;i})/2$, we have\\
\hspace{0.3in}(1)\hspace{0.1in}   $r_{ij}=Mb^2a_{ij}$,\\
\hspace{0.3in}(2)\hspace{0.1in}    $U_i=u_i+Mb_i$ is a gradient vector field,\\
\hspace{0.3in}(3)\hspace{0.1in}      $w_{ij}=s_{ij}{s^r}_j$ is written in the form\\
\hspace{1in}       $w_{ij}=b_iw_j+b_jw_i+(R-2M^2b^2)b_ib_j-b^2U_iU_j-4Ka_{ij}$,\\
 \hspace{0.3in}(4)\hspace{0.1in} $s_{ij;k}$ is written in the form \\
\hspace{1in}       $s_{ij;k}=A_{(ij)}\{a_{ik}[W_j+(R-M^2b^2)b_j]+b_iU_{jk}\}$,\\
where the symbol $A_{(ij)}$ denotes interchanges of indices $i$ and $j$ and subtraction,\\
\hspace{0.3in}(5)\hspace{0.1in} the curvature tensor $R_{hijk}$ of the associated Riemannian space has the form \\
\hspace{1in}  $^\alpha R_{hikj}=A_{(kj)}\{Ra_{hk}a_{ij}-a_{hk}U_{ij}-a_{ij}U_{hk}\}$,\\
where we put $M=a^{ij}r_{ij}/nb^2$, $u_i=b_r{s^r}_i/b^2$, $U_{ij}=U_{i;j}+U_iU_j$, $w_i=u_r{s^r}_i$, $W_i=w_i+Mb^2U_i$
and $R=4K/b^2+U_rU^r$.}\\

In [7], we characterized a Kropina space by some Riemannian space $(M, h)$ and a vector field $W$ of constant length 1
on it, and rewrite Theorem M using  $h$ and $W$. 
But, a few years after we noticed that we should replace (3) in Theorem M with the condition  
\begin{eqnarray*}
	w_{ij}=b_iw_j+b_jw_i+(R-2M^2b^2)b_ib_j-b^2U_iU_j-4Ka_{ij}+Mb^2(U_ib_j+U_jb_i) \label{}
\end{eqnarray*}
which is equivalent to the equality in Lemma 2 of [3]. Then, it follows that the condition 4 of Theorem 2 in [7]
should be omitted. Therefore, we consider the same problem in another way.

In this paper, we rewrite the condition for a Kropina space to be of constant curvature with a Riemannian metric $h$
and a vector filed $W$ and obtain the necessary and sufficient conditions for a Kropina space to be of constant 
curvature. 

 In the first section, we will show again that a Kropina space is characterized
by some Riemannian metric $h$ and a vector field $W$ of constant length 1 on it, and in the second section, we will rewrite the 
coefficients of the geodesic spray in a Kropina space by the Riemannian metric $h$ and the  vector field $W$.

In Finsler geometry, the necessary and sufficient condition for a Kropina space to be of scalar curvature is well-known 
([2]). We rewrite the condition by $h$ and $W$,  and  using the relation of $h$ and $W$ 
 we  obtain the necessary and sufficient conditions for a Kropina space to be of constant 
curvature by straightforward calculations. 
 Our main result is Theorem 4. This Theorem 4 is the corrected version of Theorem 2 in [7]. 
Therefore, this paper is the corrected version of [7].

\section{The characterization of a Kropina metric.}\label{}
 \hspace{0.2in} 
Let $(M, \alpha)$ be an $n(\ge 2)$-dimensional differential manifold endowed with a Riemannian metric $\alpha$.
A Kropina space $(M, \alpha^2/\beta)$ is a Finsler space whose fundamental function is given by $F=\alpha^2/\beta$,
 where $\alpha=\sqrt{a_{ij}(x)y^iy^j}$ and $\beta=b_i(x)y^i$.
Even though Kropina spaces can be studied in more general case ([4]),  
in this paper, we suppose that the matrix $(a_{ij})$ is 
positive definite. 

Let us remark that for a Kropina space $(M, \alpha^2/\beta)$ the Kropina metric $F=\alpha^2/\beta$ can be rewrited
 as follows:
\begin{eqnarray*}
	\frac{\alpha^2}{F^2}-\frac{\beta}{F}+\frac{b^2}{4}&=&\frac{b^2}{4}.
\end{eqnarray*}
where $b^2=a^{ij}b_ib_j$ and $(a^{ij})=(a_{ij})^{-1}$.
Let $\kappa(x)$ be a function of $(x^i)$ alone. Multiplying the above equation by $e^{\kappa(x)}$, we have
\begin{eqnarray}
     e^{\kappa(x)}a_{ij}\frac{y^i}{F}\frac{y^j}{F}
       -e^{\kappa(x)}a_{ij}\frac{y^i}{F}b^j
              +\frac{1}{4}e^{\kappa(x)}a_{ij}b^ib^j=\frac{e^{\kappa(x)}b^2}{4},
\end{eqnarray}
Defining a new Riemannian metric $h=\sqrt{h_{ij}(x)y^iy^j}$ and a vector field $W=W^i(\partial/\partial x^i)$ on $M$ by 
\begin{eqnarray}
	 h_{ij}=e^{\kappa(x)}a_{ij}\hspace{0.1in}  and \hspace{0.1in}   2W_i=e^{\kappa(x)}b_i, 
\end{eqnarray}	
where $W_i=h_{ij}W^j$, the equation (1.1) reduces to
\begin{eqnarray*}
	\bigg|\frac{y}{F}-W\bigg|=|W|.
\end{eqnarray*}
In the above equation, the symbol $| \cdot |$ denotes the length of a vector with respect to the Riemannian metric $h$.

We notice that the equation $|W|=1$ holds if and only if  the function $\kappa(x)$ satisfies the condition
\begin{eqnarray}
	e^{\kappa(x)}b^2=4.
\end{eqnarray}

Suppose that the function $\kappa(x)$ satisfies (1.3), then we have $|W|=1$ and
\begin{eqnarray}
	\bigg|\frac{y}{F}-W\bigg|=1.
\end{eqnarray}
Therefore, in each tangent space $T_xM$, the indicatrix of the Kropina metric necessarily goes through the origin.

Conversely, consider a Riemannian space $(M, h)$, where $h=\sqrt{h_{ij}(x)y^iy^j}$,
and a unit vector field $W=W^i(\partial /\partial x^i)$ on it. 
Then, we consider the metric $F$ characterized by (1.4).
Solving for $F$ from (1.4), we get
\begin{eqnarray*}
	F=\frac{|y|^2}{\{\sqrt{2}h(y,W)\}^2}.
\end{eqnarray*}

 Comparing the above equality with a Kropina metric $F=\alpha^2/\beta$, we obtain (1.2) and from the assumption
$|W|=1$ we get  (1.3).

Summarizing the above discussion, we obtain\\

{\textbf {Theorem　1}}  \hspace{0.5in}
\textit{Let  $(M, \alpha)$ be an $n(\ge 2)$-dimensional Riemannian space with a Riemannian metric
 $\alpha=\sqrt{a_{ij}(x)y^iy^j}$.
For a Kropina space $(M, F=\alpha^2/\beta)$, where  $\beta=b_i(x)y^i$, we define a new Riemannian metric $h=\sqrt{h_{ij}(x)y^iy^j}$ and a vector field 
$W=W^i(\partial/\partial x^i)$  of constant length 1 by (1.2) and (1.3). Then, the Kropina metric $F$ satisfies the equation (1.4).} 

\textit{Conversely, suppose that  $h=\sqrt{h_{ij}(x)y^iy^j}$ is  a Riemannian metric  and   $W=W^i(\partial/\partial x^i)$ is
 a   vector field of constant length 1 on $(M,h)$.
Consider the metric $F$ defined by (1.4).
Then, defining  $a_{ij}(x):=e^{-\kappa(x)}h_{ij}(x)$ and $b_i(x):=2e^{-\kappa(x)}W_i$ by (1.2) using some function
 $\kappa(x)$ of $(x^i)$ alone, we get
 $F=\alpha^2/\beta$ and it follows  the function $\kappa(x)$ satisfies (1.3).}

\section{ The coefficients of the geodesic spray in a Kropina space.  }\label{}
From the theory of Riemannian spaces, we have the following theorem:\\

{\textbf {Theorem A}}([6])\hspace{0.5in}
Let $(M, g)$ and $(M, g^*=e^{\rho}g)$, where $g=\sqrt{g_{ij}(x)y^iy^j}$ and $g^*=\sqrt{g_{ij}^*(x)y^iy^j}$ respectively,
 be two $n$-dimensional Riemannian spaces which are comformal each other.
Furthermore, let ${{\gamma_j}^i}_k$ and ${{\gamma^*_j}^i}_k$ be the coefficients of Levi-Civita connection of $(M, g)$ 
and $(M, g^*)$, respectively.
Then, we have
\begin{eqnarray*}
	g^*_{ij}=e^{2\rho}g_{ij}, \hspace{0.2in} g^{*ij}=e^{-2\rho}g^{ij},\hspace{0.2in}
	{{\gamma^*_j}^i}_k={{\gamma_j}^i}_k+\rho_j{\delta^i}_k+\rho_k{\delta^i}_j-\rho^ig_{jk},
\end{eqnarray*}
where $\rho_i=\partial \rho/\partial x^i$ and $\rho^i=g^{ij}\rho_j$.\\

From (1.2), we have $h_{ij}=e^{\kappa}a_{ij}$. Applying Theorem A , we get
\begin{eqnarray}
	{{^h \gamma_j}^i}_k={{^\alpha\gamma_j}^i}_k+\frac{1}{2}\kappa_j{\delta^i}_k+\frac{1}{2}\kappa_k{\delta^i}_j-\frac{1}{2}\kappa^ia_{jk},
\end{eqnarray}
where ${{^h \gamma_j}^i}_k$ and ${{^\alpha\gamma_j}^i}_k$ are the coefficients of Levi-Civita connection of $(M, h)$  
and $(M, \alpha)$ respectively,  $\kappa_i=\partial \kappa/\partial x^i$ and $\kappa^i=a^{ij}\kappa_j$. 
Transvecting (2.1) by $y^jy^k$, we get
\begin{eqnarray}
	{{^h \gamma_0}^i}_0={{^\alpha\gamma_0}^i}_0+\kappa_0y^i-\frac{1}{2}h_{00}\overline{\kappa}^i,
\end{eqnarray}
where $\overline{\kappa}^i=h^{ij}\kappa_j$.

We denote the covariant derivative in the Riemannian space $(M, \alpha)$ by $(_{;i})$ and put as follows:
\begin{eqnarray*}
	s_{ij}:=\frac{b_{i;j}-b_{j;i}}{2},\hspace{0.2in} r_{ij}:=\frac{b_{i;j}+b_{j;i}}{2},\hspace{0.2in}
	s_j:=b^is_{ij},\hspace{0.2in} r_j:=b^ir_{ij}.
\end{eqnarray*}

In [1],  the authors have shown  that the coefficients $G^i$ of the geodesic spray in a Finsler space
$(M, F=\alpha\phi(\beta/\alpha))$ are given by
 \begin{eqnarray}
	2G^i= {{^\alpha \gamma_0}^i}_0+2\omega \alpha {s^i}_0+2\Theta(r_{00}-2\alpha \omega s_0)
	                      \bigg(\frac{y^i}{\alpha}+\frac{\omega'}{\omega-s\omega'}b^i\bigg),
\end{eqnarray}
 where
\begin{eqnarray*}
	\omega:=\frac{\phi'}{\phi-s\phi'},\hspace{0.2in}
	\Theta:=\frac{\omega-s\omega'}{2\{1+s\omega+(b^2-s^2)\omega'\}}.
\end{eqnarray*}

For a Kropina space, we have $\phi(s)=1/s$, hence we obtain
\begin{eqnarray*}
	\phi'=-\frac{1}{s^2}, \hspace{0.2in} \omega=-\frac{1}{2s}, \hspace{0.2in} \omega'=\frac{1}{2s^2}
\end{eqnarray*}
and
\begin{eqnarray*}
	\frac{\omega'}{\omega-s\omega'}=-\frac{1}{2s},\hspace{0.1in}
	r_{00}-2\alpha \omega s_0=r_{00}+Fs_0,\hspace{0.1in}
	1+s\omega+(b^2-s^2)\omega'=\frac{b^2}{2s^2},\hspace{0.1in}
	\Theta=-\frac{s}{b^2}.
\end{eqnarray*}
Substituting the above equalities in (2.3) and using (2.2), we get
\begin{eqnarray*}
	2G^i&=&{{^h\gamma_0}^i}_0-\kappa_0y^i
	         +\frac{1}{2}h_{00}\overline{\kappa}^i-F{s^i}_0-\frac{1}{b^2}(r_{00}+Fs_0)(\frac{2}{F}y^i-b^i).
\end{eqnarray*}

From Theorem 1, for a Kropina space $(M, \alpha^2/\beta)$, a new Riemannian metric $h=\sqrt{h_{ij}(x)y^iy^j}$ and a vector 
field $W=W^i(\partial/\partial x^i)$ are defined by (1.2) and (1.3). So, the vector field $W$ satisfies the condition
$|W|=1$ and we have $F=h_{00}/2W_0$.

Therefore, we get 
\begin{eqnarray}
	2G^i={{^h\gamma_0}^i}_0+2\Phi^i,
\end{eqnarray}
where 
\begin{eqnarray}
	2\Phi^i:=-\kappa_0y^i+ \frac{1}{2}h_{00}\overline{\kappa}^i-\frac{h_{00}}{2W_0}{s^i}_0
	                               -\frac{1}{b^2}(r_{00}+\frac{h_{00}s_0}{2W_0})(\frac{4W_0}{h_{00}}y^i-b^i).
\end{eqnarray}

Using (2.1), we have
\begin{eqnarray*}
	b_{i;j} =2e^{-\kappa}W_{i||j}+e^{-\kappa}\bigg(\kappa_iW_j-\kappa_jW_i-W_r\overline{\kappa}^rh_{ij}\bigg),
\end{eqnarray*}
where the notation $(_{||i})$ stands for the  $h$-covariant derivative in the Riemannian space $(M, h)$.\\

\textbf{Remark 1}\hspace{0.7in}
\textit{We can introduce a Finsler connection $\Gamma^*=(^h{{\gamma_j}^i}_k(x),  {N_j}^i:=$
 $^h{{\gamma_j}^i}_k(x)y^k, {{C_j}^i}_k)$
associated with the linear connection $^h{{\gamma_j}^i}_k(x)$ of the Riemannian space $(M, h)$.
The \textit{$h$-covariant derivative} are defined as follows ([2]):}

\textit{For a vector field $W^i(x)$ on $M$, 
\begin{eqnarray*}
	1,\hspace{0.2in} W^i(x)_{||j}&=&\frac{\partial W^i}{\partial x^j}-\frac{\partial W^i}{\partial y^s}{N_j}^s
	               +^h{{\gamma_j}^i}_sW^s\\
	                 &=&\frac{\partial W^i}{\partial x^j}+^h{{\gamma_j}^i}_sW^s.
\end{eqnarray*}}

\textit{For a reference vector $y^i$,
\begin{eqnarray*}
	2,\hspace{0.2in} {y^i}_{||j}&=& \frac{\partial y^i}{\partial x^j}-\frac{\partial y^i}{\partial y^s}{N_j}^s
                       +^h{{\gamma_j}^i}_sy^s\\
                          &=&-{N_j}^i+{N_j}^i\\
                         &=&0.
\end{eqnarray*}}

We put 
\begin{eqnarray*}
	\texttt{R}_{ij}:=\frac{W_{i||j}+W_{j||i}}{2}, \hspace{0.1in}
	  \texttt{S}_{ij}:=\frac{W_{i||j}-W_{j||i}}{2}, \hspace{0.1in}
	&&{\texttt{R}^i}_j:=h^{ir}\texttt{R}_{rj},\hspace{0.1in}
        {\texttt{S}^i}_j:=h^{ir}\texttt{S}_{rj}\\	
	\texttt{R}_{i}:=W^r\texttt{R}_{ri},    \hspace{0.1in}
	 \texttt{S}_{i}:=W^r\texttt{S}_{ri},    \hspace{0.1in}
	&&\texttt{R}^{i}:=h^{ir}\texttt{R}_{r}, \hspace{0.1in}
       \texttt{S}^{i}:=h^{ir}\texttt{S}_{r}.
\end{eqnarray*}

It follows
\begin{eqnarray*}
	r_{ij}=2e^{-\kappa} \bigg(\texttt{R}_{ij}-\frac{1}{2}W_r\overline{\kappa}^rh_{ij}\bigg),\hspace{0.1in}
	s_{ij}=2e^{-\kappa} \bigg(\texttt{S}_{ij}+\frac{\kappa_iW_j-\kappa_jW_i}{2}\bigg).
\end{eqnarray*}
Furthermore, we get
\begin{eqnarray*}
	{s^i}_j=2{\texttt{S}^i}_j+\overline{\kappa}^iW_j-\kappa_jW^i,&&\hspace{0.1in}
	{s^i}_0=2{\texttt{S}^i}_0+W_0\overline{\kappa}^i-\kappa_0W^i\\
	s_i=2e^{-\kappa} \bigg(2\texttt{S}_i+W_r\overline{\kappa}^rW_i-\kappa_i\bigg),&&\hspace{0.1in}
	s_0=2e^{-\kappa} \bigg(2\texttt{S}_0+W_r\overline{\kappa}^rW_0-\kappa_0\bigg) ,\\
	r_{00}=2e^{-\kappa} \bigg(\texttt{R}_{00}-\frac{1}{2}W_r\overline{\kappa}^rh_{00}\bigg) &&\hspace{0.1in} 
	b^i=a^{ir}b_r=e^\kappa h^{ir}\frac{2W_r}{e^\kappa}=2W^i.	
\end{eqnarray*}
Substituting the above equality in (2.5), we have
\begin{eqnarray}
	2\Phi^i=\frac{h_{00}}{W_0}(\texttt{S}_0W^i-{\texttt{S}^i}_0)
                       +(\texttt{R}_{00}W^i-2\texttt{S}_0y^i)
                      -\frac{2W_0}{h_{00}}\texttt{R}_{00}y^i.	
\end{eqnarray}
Multiplying now the above equality by $2h_{00}W_0$, we get
\begin{eqnarray*}
	&&4h_{00}W_0 \Phi^i
 =2(h_{00})^2(\texttt{S}_0W^i-{\texttt{S}^i}_0)+2h_{00}W_0(\texttt{R}_{00}W^i-2\texttt{S}_0y^i)
	             -4(W_0)^2\texttt{R}_{00}y^i
\end{eqnarray*}
and by putting
\begin{eqnarray*}
	A_{(1)}^i:=2(\texttt{S}_0W^i-{\texttt{S}^i}_0),\hspace{0.1in}
	A_{(2)}^i:=2(\texttt{R}_{00}W^i-2\texttt{S}_0y^i),\hspace{0.1in}
	A_{(3)}^i:=-4\texttt{R}_{00}y^i,	
\end{eqnarray*}
it follows 
\begin{eqnarray}
	4h_{00}W_0 \Phi^i	=(h_{00})^2A_{(1)}^i +h_{00}W_0A_{(2)}^i        +(W_0)^2A_{(3)}^i.
\end{eqnarray}

\section{The necessary and sufficient conditions for a Kropina space to be of constant curvature.}\label{}
	
\hspace{0.2in}
In this section, we consider a Kropina space $(M, F=\alpha^2/\beta)$ of constant curvatue $K$, where 
$\alpha=\sqrt{a_{ij}(x)y^iy^j}$ and $\beta=b_i(x)y^i$. Furthermore, we suppose that  the matrix $(a_{ij})$ is always
 positive definite and that the dimension $n$ is more than or equal two.
Hence, it follows  that $\alpha^2$ is not divisible by $\beta$.  This is an important relation and 
is equivalent to that $h_{00}$ is not divisible by $W_0$.
Using these, we will obtain  the necessary and sufficient conditions for a Kropina space to be of 
constant curvature. 

\subsection{The curvature tensor of a  Kropina space.}\label{}
\hspace{0.2in}
Let ${{R_j}^i}_{kl}$ be the $h$-curvature tensor of Cartan connection in Finsler space. 
The Berwald's spray curvature tensor is
\begin{eqnarray}
	^{(b)}{{R_j}^i}_{kl}&=&A_{(kl)}\bigg\{\frac{\partial {{G_j}^i}_k}{\partial x^l}+{{G_j}^r}_k{{G_r}^i}_l\bigg\}.	           
\end{eqnarray}
It is well-known that the equality ${{R_0}^i}_{kl}=^{(b)}{{R_0}^i}_{kl}$ holds good ([2]).

From $2G^i={{^h\gamma_0}^i}_0+2\Phi^i$, it follows
\begin{eqnarray*}
	{G^i}_j={{^h\gamma_0}^i}_j+{\Phi^i}_j \hspace{0.2in} and \hspace{0.2in}
	{{G_j}^i}_k={{^h\gamma_j}^i}_k+{{\Phi_j}^i}_k.
\end{eqnarray*}
Substituting the above equalities in (3.1), we get
\begin{eqnarray*}
	{{R_j}^i}_{kl}={{^hR_j}^i}_{kl}+ A_{(kl)}\bigg\{ {{\Phi_j}^i}_{k||l}+{{\Phi_j}^r}_k{{\Phi_r}^i}_l\bigg\}.
\end{eqnarray*}

The following Proposition is well-known ([2]):\\

\textbf{Proposition 1}\hspace{0.5in}
  \textit{The necessary and sufficient condition for a Finsler space $(M, F)$ to be of scalar curvature $K$ is that the equality 
\begin{eqnarray}
	 {{R_0}^i}_{0l}=KF^2 ({\delta^i}_l-l^il_l),
\end{eqnarray}
where $l^i=y^i/F$ and $l_l=\partial F/\partial y^l$,  holds.} \\

If the equality (3.2) holds good and $K$ is constant, then the Finsler space is called of {\it constant curvature} $K$.

For a Kropina space of constant curvatue $K$, the equality (3.2) holds and $K$ is constant.

Since $F=\alpha^2/\beta$ is written as $F=h_{00}/(2W_0)$, we have
\begin{eqnarray*}
	{\delta^i}_l-l^il_l={\delta^i}_l-\frac{2W_0h_{0l}-h_{00}W_l}{h_{00}W_0}y^i.
\end{eqnarray*}

Using the curvature we obtained above, we have
\begin{eqnarray*}
	{{R_0}^i}_{0l}={{^hR_0}^i}_{0l}+  2{\Phi^i}_{||l}-  {\Phi^i}_{l||0}
	                                 +2\Phi^r{{\Phi_r}^i}_l -{\Phi^r}_l{\Phi^i}_r.
\end{eqnarray*}
Substituting the above equalities in (3.2), we get
\begin{eqnarray}
	KF^2\bigg({\delta^i}_l-\frac{2W_0h_{0l}-h_{00}W_l}{h_{00}W_0}y^i\bigg)={{^hR_0}^i}_{0l}+  2{\Phi^i}_{||l}-  {\Phi^i}_{l||0}
	                                 +2\Phi^r{{\Phi_r}^i}_l -{\Phi^r}_l{\Phi^i}_r.
\end{eqnarray}

\subsection{Rewriting the equation (3.3) using $h_{00}$ and $W_0$. }\label{}
\textbf{(1)The calculations for ${\Phi^i}_{||l}$.}

First, applying the $h$-covariant derivative  $_{||l}$ to (2.7), it follows 
\begin{eqnarray*}
	&&4h_{00}W_{0||l} \Phi^i+4h_{00}W_0 {\Phi^i}_{||l}\\
	&=&(h_{00})^2{A_{(1)}^i}_{||l} +h_{00}W_0{A_{(2)}^i}_{||l} 	                                  
	  +h_{00}W_{0||l}A_{(2)}^i+(W_0)^2{A_{(3)}^i}_{||l} 	+2W_0W_{0||l}A_{(3)}^i  \nonumber	
\end{eqnarray*}
and using again  (2.7), we get
\begin{eqnarray*}
	&&4h_{00}(W_0)^2 {\Phi^i}_{||l}\\
	&=&(h_{00})^2W_0{A_{(1)}^i}_{||l} 
	     -(h_{00})^2A_{(1)}^iW_{0||l}
	     +h_{00}(W_0)^2{A_{(2)}^i}_{||l} \\
	 &&\hspace{2in}        +(W_0)^3{A_{(3)}^i}_{||l}
	               +(W_0)^2A_{(3)}^iW_{0||l}.	
\end{eqnarray*}
Putting now
\begin{eqnarray*}	
	&&B_{(1)l}^i={A_{(1)}^i}_{||l} ,\hspace{0.1in}
      B_{(21)l}^i=-A_{(1)}^iW_{0||l}, \hspace{0.1in} 
       B_{(22)l}^i={A_{(2)}^i}_{||l} ,\hspace{0.1in}
        B_{(3)l}^i={A_{(3)}^i}_{||l},\hspace{0.1in}
                         B_{(4)l}^i=A_{(3)}^iW_{0||l},
\end{eqnarray*}
we have
\begin{eqnarray}
	&&4h_{00}(W_0)^2 {\Phi^i}_{||l}
	=(h_{00})^2W_0B_{(1)l}^i
	   +(h_{00})^2B_{(21)l}^i      \\
       &&\hspace{2in} +h_{00}(W_0)^2B_{(22)l}^i    +(W_0)^3B_{(3)l}^i
                     +(W_0)^2B_{(4)l}^i.\nonumber
\end{eqnarray}

\textbf{(2)The calculations for ${\Phi^i}_l$.}

Secondly, derivating (2.7) by $y^l$, we get
\begin{eqnarray*}
	&&8h_{0l}W_0 \Phi^i+4h_{00}W_l \Phi^i+4h_{00}W_0 {\Phi^i}_l\\
	&=& (h_{00})^2A_{(1)}^i._l      +h_{00}W_0A_{(2).l}^i 
	             +h_{00}(4h_{0l}A_{(1)}^i+W_lA_{(2)}^i)\\
	&&\hspace{2in}   +(W_0)^2A_{(3)}^i._l   +2W_0(W_lA_{(3)}^i+h_{0l}A_{(2)}^i),
\end{eqnarray*}
where the notation  $(_{.l})$ stands for the derivation by $y^l$.

Using the above equality and (2.7), we obtain
\begin{eqnarray}
	 &&4(h_{00})^2(W_0)^2 {\Phi^i}_l\\
	      &=&(h_{00})^3W_0C_{(0)l}^i
            +(h_{00})^3C_{(11)l}^i
            +(h_{00})^2(W_0)^2C_{(12)l}^i       \nonumber   \\
      &&\hspace{1in}  +(h_{00})^2W_0C_{(21)l}^i
               +h_{00}(W_0)^3C_{(22)l}^i
                +h_{00}(W_0)^2C_{(3)l}^i
	       +(W_0)^3C_{(4)l}^i,\nonumber
\end{eqnarray}
where
\begin{eqnarray*}
	&&C_{(0)l}^i:=A_{(1)}^i._l=2(\texttt{S}_lW^i-{\texttt{S}^i}_l), \hspace{0.1in}
      C_{(11)l}^i:=-A_{(1)}^iW_l=-2(\texttt{S}_0W^i-{\texttt{S}^i}_0)W_l,\\
      &&C_{(12)l}^i:=A_{(2).l}^i
                 =4( \texttt{R}_{0l}W^i-\texttt{S}_ly^i-\texttt{S}_0{\delta^i}_l),\hspace{0.1in}  
      C_{(21)l}^i:=2A_{(1)}^ih_{0l}=4(\texttt{S}_0W^i-{\texttt{S}^i}_0)h_{0l},\\
      &&C_{(22)l}^i:=A_{(3)}^i._l=-4(2\texttt{R}_{0l}y^i+\texttt{R}_{00}{\delta^i}_l),\hspace{0.1in} 
      C_{(3)l}^i:=W_lA_{(3)}^i=-4\texttt{R}_{00}W_ly^i,\\
      &&C_{(4)l}^i:=-2A_{(3)}^ih_{0l}=8\texttt{R}_{00}h_{0l}y^i.
\end{eqnarray*}

\textbf{(3)The calculations for ${\Phi^i}_{l||0}$.}

Applying the $h$-covariant derivaive ${_{||0}}$ to (3.5), we get
\begin{eqnarray*}
	 &&8(h_{00})^2W_0 W_{0||0}{\Phi^i}_l+4(h_{00})^2(W_0)^2 {\Phi^i}_{l||0}\\
	&=&(h_{00})^3W_0C_{(0)l||0}^i
	   +(h_{00})^3(W_{0||0}C_{(0)l}^i+C_{(11)l||0}^i)       \\
        && +(h_{00})^2(W_0)^2C_{(12)l||0}^i
      +(h_{00})^2W_0(2W_{0||0}C_{(12)l}^i+C_{(21)l||0}^i)\\
      &&+(h_{00})^2W_{0||0}C_{(21)l}^i 
         +h_{00}(W_0)^3C_{(22)l||0}^i 
          +h_{00}(W_0)^2(3W_{0||0}C_{(22)l}^i+C_{(3)l||0}^i)\\
       &&+2h_{00}W_0W_{0||0}C_{(3)l}^i 
               +(W_0)^3C_{(4)l||0}^i 
	        +3(W_0)^2W_{0||0}C_{(4)l}^i.	
\end{eqnarray*}

Using the above equality and (3.5), we obtain
\begin{eqnarray}
	&&4(h_{00})^2(W_0)^3 {\Phi^i}_{l||0}\\
	&=&(h_{00})^3(W_0)^2D_{(1)l}^i
	               +(h_{00})^3W_0D_{(21)l}^i
	                     +(h_{00})^3D_{(31)l}^i\nonumber\\
	&&\hspace{0.5in}   +(h_{00})^2(W_0)^3D_{(22)l}^i
	              +(h_{00})^2(W_0)^2D_{(32)l}^i
	                  +(h_{00})^2W_0D_{(41)l}^i\nonumber\\
	&&\hspace{0.5in}   +h_{00}(W_0)^4D_{(33)l}^i
                       +h_{00}(W_0)^3 D_{(42)l}^i
                         +(W_0)^4D_{(5)l}^i
                         +(W_0)^3D_{(6)l}^i,\nonumber
\end{eqnarray}
where
\begin{eqnarray*}
	&&D_{(1)l}^i=C_{(0)l||0}^i,\hspace{0.2in}
	        D_{(21)l}^i=C_{(11)l||0}^i-W_{0||0}C_{(0)l}^i,\hspace{0.2in}
               D_{(31)l}^i=-2W_{0||0}C_{(11)l}^i,   \\
	&&   D_{(22)l}^i=C_{(12)l||0}^i,\hspace{0.2in}
	           D_{(32)l}^i=C_{(21)l||0}^i , \hspace{0.2in}
	             D_{(41)l}^i=-W_{0|0}C_{(21)l}^i, \\
	&&D_{(33)l}^i=C_{(22)l||0}^i , \hspace{0.2in}
                 D_{(42)l}^i=W_{0||0}C_{(22)l}^i+C_{(3)l||0}^i, \hspace{0.2in}
                     D_{(5)l}^i= C_{(4)l||0}^i  , \hspace{0.2in}
                   D_{(6)l}^i=    W_{0||0}C_{(4)l}^i.
\end{eqnarray*}

\textbf{(4)The calculations for ${\Phi^r}_l{\Phi^i}_r$.}

Using (3.5), we have
\begin{eqnarray}
	 &&16(h_{00})^4(W_0)^4 {\Phi^r}_l{\Phi^i}_r\\
 &=&(h_{00})^6(W_0)^2E_{(01)l}^i
      +(h_{00})^6W_0E_{(11)l}^i  
          +(h_{00})^6E_{(21)l}^i
               +(h_{00})^5(W_0)^3E_{(12)l}^i  \nonumber\\
   &&+(h_{00})^5(W_0)^2E_{(22)l}^i
          +(h_{00})^5W_0E_{(31)l}^i 
              +(h_{00})^4(W_0)^4E_{(23)l}^i
                   +(h_{00})^4(W_0)^3E_{(32)l}^i \nonumber\\
    &&+(h_{00})^4(W_0)^2E_{(41)l}^i 
              +(h_{00})^3(W_0)^5E_{(33)l}^i          
                 +(h_{00})^3(W_0)^4E_{(42)l}^i
                      +(h_{00})^3(W_0)^3E_{(51)l}^i  \nonumber\\
     &&+(h_{00})^2(W_0)^6E_{(43)l}^i
          +(h_{00})^2(W_0)^5E_{(52)l}^i  
                +(h_{00})^2(W_0)^4E_{(61)l}^i
                       +h_{00}(W_0)^6E_{(62)l}^i  \nonumber\\
     &&+h_{00}(W_0)^5E_{(7)l}^i
                +(W_0)^6E_{(8)l}^i,  \nonumber
\end{eqnarray}
where
\begin{eqnarray*}
E_{(0)l}^i:&=&C_{(0)r}^iC_{(0)l}^r,\hspace{0.1in}            
                  E_{(11)l}^i:=C_{(11)r}^iC_{(0)l}^r+C_{(0)r}^iC_{(11)l}^r,\hspace{0.1in}
                  E_{(21)l}^i:=C_{(11)r}^iC_{(11)l}^r,\\           
E_{(12)l}^i:&=&C_{(0)r}^iC_{(12)l}^r  +C_{(12)r}^iC_{(0)l}^r ,\\
E_{(22)l}^i:&=&C_{(12)r}^iC_{(11)l}^r    +C_{(11)r}^iC_{(12)l}^r+C_{(21)r}^iC_{(0)l}^r +C_{(0)r}^iC_{(21)l}^r ,\\
E_{(31)l}^i:&=&C_{(21)r}^iC_{(11)l}^r +C_{(11)r}^iC_{(21)l}^r, \hspace{0.1in}
                 E_{(23)l}^i:=C_{(12)r}^iC_{(12)l}^r+C_{(22)r}^iC_{(0)l}^r       +C_{(0)r}^iC_{(22)l}^r,   \\
E_{(32)l}^i:&=&C_{(21)r}^iC_{(12)l}^r +C_{(12)r}^iC_{(21)l}^r+C_{(3)r}^iC_{(0)l}^r+C_{(22)r}^iC_{(11)l}^r
                      +C_{(11)r}^iC_{(22)l}^r     +C_{(0)r}^iC_{(3)l}^r \\
E_{(41)l}^i:&=&C_{(3)r}^iC_{(11)l}^r+C_{(11)r}^iC_{(3)l}^r  +C_{(21)r}^iC_{(21)l}^r, \hspace{0.1in}
E_{(33)l}^i:=C_{(22)r}^iC_{(12)l}^r+C_{(12)r}^iC_{(22)l}^r , \\
E_{(42)l}^i:&=&C_{(4)r}^iC_{(0)l}^r+C_{(3)r}^iC_{(12)l}^r+C_{(22)r}^iC_{(21)l}^r
               +C_{(21)r}^iC_{(22)l}^r+C_{(12)r}^iC_{(3)l}^r +C_{(0)r}^iC_{(4)l}^r,\\
E_{(51)l}^i:&=&C_{(3)r}^iC_{(21)l}^r+C_{(21)r}^iC_{(3)l}^r+C_{(4)r}^iC_{(11)l}^r +C_{(11)r}^iC_{(4)l}^r,\hspace{0.1in}
 E_{(43)l}^i:=C_{(22)r}^iC_{(22)l}^r , \\
 E_{(52)l}^i:&=&C_{(4)r}^iC_{(12)l}^r+C_{(12)r}^iC_{(4)l}^r
                +C_{(3)r}^iC_{(22)l}^r  +C_{(22)r}^iC_{(3)l}^r,   \\
 E_{(61)l}^i:&=&C_{(3)r}^iC_{(3)l}^r +C_{(4)r}^iC_{(21)l}^r  +C_{(21)r}^iC_{(4)l}^r,\hspace{0.1in}
 E_{(62)l}^i:=C_{(4)r}^iC_{(22)l}^r  +C_{(22)r}^iC_{(4)l}^r,\\
 E_{(7)l}^i:&=&C_{(4)r}^iC_{(3)l}^r+C_{(3)r}^iC_{(4)l}^r,\hspace{0.1in}             
 E_{(8)l}^i:=C_{(4)r}^iC_{(4)l}^r.
\end{eqnarray*}

\textbf{(5)The calculations for ${\Phi^r}{{\Phi_r}^i}_l$.}

Derivating (3.5) by $y^r$, we get
\begin{eqnarray*}
	 &&16h_{00}h_{0r}(W_0)^2 {\Phi^i}_l+8(h_{00})^2W_0W_r {\Phi^i}_l+4(h_{00})^2(W_0)^2 {{\Phi_r}^i}_l\\
	&=&6(h_{00})^2h_{0r}W_0C_{(0)l}^i+(h_{00})^3W_rC_{(0)l}^i+(h_{00})^3W_0C_{(0)l}^i._r\\
            &&+6(h_{00})^2h_{0r}C_{(11)l}^i+(h_{00})^3C_{(11)l}^i._r\\
            &&+4h_{00}h_{0r}(W_0)^2C_{(12)l}^i+2(h_{00})^2W_0W_rC_{(12)l}^i+(h_{00})^2(W_0)^2C_{(12)l}^i._r\\
            &&+4h_{00}h_{0r}W_0C_{(21)l}^i+(h_{00})^2W_rC_{(21)l}^i+(h_{00})^2W_0C_{(21)l}^i._r\\
            && +2h_{0r}(W_0)^3C_{(22)l}^i+3h_{00}(W_0)^2W_rC_{(22)l}^i+h_{00}(W_0)^3C_{(22)l}^i._r\\
            &&+2h_{0r}(W_0)^2C_{(3)l}^i+2h_{00}W_0W_rC_{(3)l}^i+h_{00}(W_0)^2C_{(3)l}^i._r\\
	      &&+3(W_0)^2W_rC_{(4)l}^i+(W_0)^3C_{(4)l}^i._r.
\end{eqnarray*}
Using the above equality, (3.5) and  $C_{(0)l}^i._r=0$, we have
\begin{eqnarray*}
	&&4(h_{00})^3(W_0)^3 {{\Phi_r}^i}_l\\
&=&(h_{00})^4W_0H_{(01)rl}^i
            +(h_{00})^4H_{(11)rl}^i
            +(h_{00})^3(W_0)^3H_{(02)rl}^i
           +(h_{00})^3(W_0)^2H_{(12)rl}^i \\
&&+(h_{00})^3W_0H_{(21)rl}^i
            +(h_{00})^2(W_0)^4H_{(13)rl}^i
            +(h_{00})^2(W_0)^3H_{(22)rl}^i\\
&& +h_{00}(W_0)^4H_{(3)rl}^i
            +h_{00}(W_0)^3H_{(4)rl}^i 
            +(W_0)^4H_{(5)rl}^i,
\end{eqnarray*}
where
\begin{eqnarray*}
&&H_{(01)rl}^i=C_{(11)l}^i._r -W_rC_{(0)l}^i, \hspace{0.1in}
               H_{(11)rl}^i=-2W_rC_{(11)l}^i,\hspace{0.1in}
              H_{(02)rl}^i=C_{(12)l.r}^i,\\
&&H_{(12)rl}^i=2h_{0r}C_{(0)l}^i+C_{(21)l}^i._r,  \hspace{0.1in}
               H_{(21)rl}^i=-W_rC_{(21)l}^i+2h_{0r}C_{(11)l}^i,   \\
&&H_{(13)rl}^i=C_{(22)l}^i._r, \hspace{0.1in}
             H_{(22)rl}^i=W_rC_{(22)l}^i+C_{(3)l}^i._r, \hspace{0.1in}
              H_{(3)rl}^i=C_{(4)l}^i._r-2h_{0r}C_{(22)l}^i,\\
&&H_{(4)rl}^i=W_rC_{(4)l}^i-2h_{0r}C_{(3)l}^i,\hspace{0.1in}
             H_{(5)rl}^i=-4h_{0r}C_{(4)l}^i.	
\end{eqnarray*}

Using the above equality and (2.7), we get
\begin{eqnarray}
	&&16(h_{00})^4(W_0)^4  \Phi^r{{\Phi_r}^i}_l\\
&=&(h_{00})^6W_0J_{(11)l}^i 
        +(h_{00})^6J_{(21)l}^i
        +(h_{00})^5(W_0)^3J_{(12)l}^i 
          +(h_{00})^5(W_0)^2J_{(22)l}^i\nonumber\\
    &&+(h_{00})^5W_0J_{(31)l}^i
          +(h_{00})^4(W_0)^4J_{(23)l}^i
           +(h_{00})^4(W_0)^3J_{(32)l}^i
            +(h_{00})^4(W_0)^2J_{(41)l}^i  \nonumber\\
     &&+(h_{00})^3(W_0)^5J_{(33)l}^i
             + (h_{00})^3(W_0)^4J_{(42)l}^i
             + (h_{00})^3(W_0)^3J_{(51)l}^i  
              +(h_{00})^2(W_0)^6J_{(43)l}^i\nonumber\\
     &&+(h_{00})^2(W_0)^5J_{(52)l}^i 
              + (h_{00})^2(W_0)^4J_{(61)l}^i
               +h_{00}(W_0)^6J_{(62)l}^i  
               +h_{00}(W_0)^5J_{(71)l}^i
               +(W_0)^6J_{(8)l}^i, \nonumber
\end{eqnarray}
where
\begin{eqnarray*}
&&J_{(11)l}^i=A_{(1)}^rH_{(01)rl}^i,\hspace{0.1in}
        J_{(21)l}^i=A_{(1)}^rH_{(11)rl}^i,\hspace{0.1in}
              J_{(12)l}^i=A_{(1)}^rH_{(02)rl}^i,\\
&&J_{(22)l}^i=A_{(1)}^rH_{(12)rl}^i+A_{(2)}^rH_{(01)rl}^i,\hspace{0.1in}
              J_{(31)l}^i=A_{(1)}^rH_{(21)rl}^i+A_{(2)}^rH_{(11)rl}^i,\\
&&J_{(23)l}^i=A_{(1)}^rH_{(13)rl}^i+A_{(2)}^rH_{(02)rl}^i,\hspace{0.1in}
         J_{(32)l}^i=A_{(1)}^rH_{(22)rl}^i +A_{(2)}^rH_{(12)rl}^i+A_{(3)}^rH_{(01)rl}^i,\\
&&J_{(41)l}^i=A_{(2)}^rH_{(21)rl}^i+A_{(3)}^rH_{(11)rl}^i,\hspace{0.1in}
                J_{(33)l}^i=A_{(2)}^rH_{(13)rl}^i+A_{(3)}^rH_{(02)rl}^i,\\
&&J_{(42)l}^i=A_{(1)}^rH_{(3)rl}^i
          +A_{(2)}^rH_{(22)rl}^i+A_{(3)}^rH_{(12)rl}^i,\\
&&J_{(51)l}^i=A_{(1)}^rH_{(4)rl}^i+A_{(3)}^rH_{(21)rl}^i,\hspace{0.1in}
                    J_{(43)l}^i=A_{(3)}^rH_{(13)rl}^i,\\
&&J_{(52)l}^i=A_{(2)}^rH_{(3)rl}^i+A_{(3)}^rH_{(22)rl}^i,\hspace{0.1in}
              J_{(61)l}^i=A_{(2)}^rH_{(4)rl}^i+A_{(1)}^rH_{(5)rl}^i,\\
&&J_{(62)l}^i=A_{(3)}^rH_{(3)rl}^i,\hspace{0.1in}
             J_{(7)l}^i=A_{(2)}^rH_{(5)rl}^i+ A_{(3)}^rH_{(4)rl}^i,\hspace{0.1in}
   J_{(8)l}^i=A_{(3)}^rH_{(5)rl}^i.
\end{eqnarray*}

\textbf{(6)The main relation.}

Multiplying (3.3) by $16(h_{00})^4(W_0)^4$ and  using $F^2=(h_{00})^2/\{4(W_0)^2\}$, we have the equality
\begin{eqnarray*}
	&&4K(h_{00})^6(W_0)^2 {h^i}_l
	= 16(h_{00})^4(W_0)^4 \cdot {{^hR_0}^i}_{0l}+8(h_{00})^3(W_0)^2  \cdot 4h_{00}(W_0)^2{\Phi^i}_{||l}    \\
	       &&- 4(h_{00})^2W_0\cdot 4(h_{00})^2(W_0)^3{\Phi^i}_{l||0}
	           +2 \cdot 16(h_{00})^4(W_0)^4\Phi^r{{\Phi_r}^i}_l -16(h_{00})^4(W_0)^4{\Phi^r}_l{\Phi^i}_r.
\end{eqnarray*}
Substituting (3.4), (3.5), (3.6), (3.7) and (3.8) in the above equality, we get
\begin{eqnarray*}
&&4K(h_{00})^6(W_0)^2{\delta^i}_l-8K(h_{00})^5(W_0)^2h_{0l}y^i+4K(h_{00})^6W_0W_ly^i\\	
&=&-(h_{00})^6(W_0)^2E_{(0)l}^i
            +(h_{00})^6W_0(2J_{(11)l}^i -E_{(11)l}^i)\\
&&+(h_{00})^6(2J_{(21)l}^i-E_{(21)l}^i )
            +(h_{00})^5(W_0)^3(8B_{(1)l}^i-4D_{(1)l}^i+2J_{(12)l}^i-E_{(12)l}^i) \\
&&+(h_{00})^5(W_0)^2(-4D_{(21)l}^i+8B_{(21)l}^i+2J_{(22)l}^i-E_{(22)l}^i)\\
&&+(h_{00})^5W_0(2J_{(31)l}^i -4D_{(31)l}^i-E_{(31)l}^i)\\
&&+(h_{00})^4(W_0)^4(16 {{^hR_0}^i}_{0l}+8B_{(22)l}^i-4D_{(22)l}^i+2J_{(23)l}^i-E_{(23)l}^i)\\
&& +(h_{00})^4(W_0)^3(2J_{(32)l}^i-4D_{(32)l}^i -E_{(32)l}^i )
      +(h_{00})^4(W_0)^2(2J_{(41)l}^i -E_{(41)l}^i-4D_{(41)l}^i)\\
&&+(h_{00})^3(W_0)^5(8B_{(3)l}^i-4D_{(33)l}^i+2J_{(33)l}^i-E_{(33)l}^i) \\
&&+ (h_{00})^3(W_0)^4(2J_{(42)l}^i+8B_{(4)l}^i-4 D_{(42)l}^i-E_{(42)l}^i)\\
&&+ (h_{00})^3(W_0)^3(2J_{(51)l}^i-E_{(51)l}^i)
        +(h_{00})^2(W_0)^6(2J_{(43)l}^i-E_{(43)l}^i)\\
&&+(h_{00})^2(W_0)^5(2J_{(52)l}^i-4D_{(5)l}^i-E_{(52)l}^i)  
          + (h_{00})^2(W_0)^4(2J_{(61)l}^i-4D_{(6)l}^i-E_{(61)l}^i)\\
&&+h_{00}(W_0)^6(2J_{(62)l}^i-E_{(62)l}^i )
         +h_{00}(W_0)^5(2J_{(71)l}^i-E_{(7)l}^i)
               +(W_0)^6(2J_{(8)l}^i-E_{(8)l}^i).
\end{eqnarray*}

Taking into account the equalities $J_{(21)l}^i=0$ and $E_{(21)l}^i=0$, we have
\begin{eqnarray}
             (h_{00})^4{P_{(5)}^i}_l+(h_{00})^2{Q_{(9)}^i}_l+(W_0)^4{R_{(9)}^i}_l=0,
\end{eqnarray}
where
\begin{eqnarray*}
	{P_{(5)}^i}_l&=&(h_{00})^2W_0(-E_{(0)l}^i-4K{\delta^i}_l)\\
                     &&\hspace{0.2in}+(h_{00})^2(2J_{(11)l}^i -E_{(11)l}^i-4KW_ly^i)\nonumber\\
                     &&+h_{00}(W_0)^2(8B_{(1)l}^i-4D_{(1)l}^i+2J_{(12)l}^i-E_{(12)l}^i) \nonumber\\
                     &&\hspace{0.2in}+h_{00}W_0(-4D_{(21)l}^i+8B_{(21)l}^i+2J_{(22)l}^i-E_{(22)l}^i+8Kh_{0l}y^i)\nonumber\\  
                     &&+(W_0)^3(16 {{^hR_0}^i}_{0l}+8B_{(22)l}^i-4D_{(22)l}^i+2J_{(23)l}^i-E_{(23)l}^i)\nonumber\\
                     &&\hspace{0.2in} +(W_0)^2(2J_{(32)l}^i-4D_{(32)l}^i -E_{(32)l}^i ),\nonumber\\
     {Q_{(9)}^i}_l&=&(h_{00})^3(2J_{(31)l}^i -4D_{(31)l}^i-E_{(31)l}^i)\\
                     &&\hspace{0.4in}+(h_{00})^2W_0(2J_{(41)l}^i -E_{(41)l}^i-4D_{(41)l}^i)\nonumber\\
                      &&+h_{00}(W_0)^4(8B_{(3)l}^i-4D_{(33)l}^i+2J_{(33)l}^i-E_{(33)l}^i) \nonumber\\
                      &&\hspace{0.2in}+ h_{00}(W_0)^3(2J_{(42)l}^i+8B_{(4)l}^i-4 D_{(42)l}^i-E_{(42)l}^i)\nonumber\\
                     &&\hspace{0.4in}+ h_{00}(W_0)^2(2J_{(51)l}^i-E_{(51)l}^i)\nonumber\\
                      &&+(W_0)^5(2J_{(43)l}^i-E_{(43)l}^i)\nonumber\\
                      &&\hspace{0.2in}+(W_0)^4(2J_{(52)l}^i-4D_{(5)l}^i-E_{(52)l}^i)  \nonumber\\
                       &&\hspace{0.4in}+ (W_0)^3(2J_{(61)l}^i-4D_{(6)l}^i-E_{(61)l}^i)\nonumber\\
      {R_{(9)}^i}_l&=&h_{00}W_0(2J_{(62)l}^i-E_{(62)l}^i )\\
                        &&\hspace{0.6in}+h_{00}(2J_{(7)l}^i-E_{(7)l}^i)\nonumber\\
                             &&+W_0(2J_{(8)l}^i-E_{(8)l}^i) . \nonumber
\end{eqnarray*}
We call ${P_{(5)}^i}_l$, ${Q_{(9)}^i}_l$ and  ${R_{(9)}^i}_l$ the  \textit{curvature part}, 
the \textit{vanishing part} and the \textit{Killing part}, respectively.\\

\textbf{Proposition 2}\hspace{0.5in}
\textit{The necessary and sufficient condition for a Kropina space $(M, F=\alpha^2/\beta=h_{00}/2W_0)$ to be of constant curvature
$K$ is that (3.9) holds good.}

\subsection{The Killing part.}\label{}
We consider the Killing part  ${R_{(9)}^i}_l$ and obtain the conclusion that the vector field $W$ is Killing. 

First, we have	
\begin{eqnarray*}
	J_{(62)l}^i &=&-64\texttt{R}_{00}(2\texttt{R}_{00}h_{0l}+h_{00}\texttt{R}_{0l})y^i
	                                 -32(\texttt{R}_{00})^2h_{00}{\delta^i}_l,\\	
	E_{(62)l}^i&=&-64\texttt{R}_{00} h_{00}\texttt{R}_{0l}y^i
	                 -128 (\texttt{R}_{00})^2h_{0l}y^i,\\
	J_{(7)l}^i &=&-32\texttt{R}_{00}\{(3\texttt{R}_{00}W_0-4\texttt{S}_0h_{00}) h_{0l}
                                      + h_{00}\texttt{R}_{00}W_l\}y^i,\\
	E_{(7)l}^i &=&-32(\texttt{R}_{00})^2(h_{00}W_ly^i+W_0h_{0l}y^i),\\
	2J_{(8)l}^i-E_{(8)l}^i&=&192(\texttt{R}_{00})^2h_{00}h_{0l}y^i.\nonumber
\end{eqnarray*}
Using the above equalities, we get
\begin{eqnarray*}
	{R_{(9)}^i}_l  =-32h_{00}\texttt{R}_{00}\bigg(W_0(2\texttt{R}_{00}h_{00}{\delta^i}_l
	                                 + 2h_{00}\texttt{R}_{0l}y^i
	                               +7\texttt{R}_{00}h_{0l}y^i)
	                        -8\texttt{S}_0h_{00} h_{0l}  y^i
                                               +\texttt{R}_{00}h_{00}W_ly^i\bigg).
\end{eqnarray*}
Substituting the above equality in (3.9) and dividing it by $W_0h_{00}$, we get
\begin{eqnarray}
	 &&(h_{00})^3{P_{(5)}^i}_l+h_{00}{Q_{(9)}^i}_l\\
	 &&\hspace{0.5in} -32(W_0)^4\texttt{R}_{00}\bigg(W_0(2\texttt{R}_{00}h_{00}{\delta^i}_l
	                                 + 2h_{00}\texttt{R}_{0l}y^i
	                               +7\texttt{R}_{00}h_{0l}y^i)  \nonumber\\
	       &&\hspace{2.5in}                 -8\texttt{S}_0h_{00} h_{0l}  y^i
                                               +\texttt{R}_{00}h_{00}W_ly^i\bigg)=0.\nonumber
\end{eqnarray}\\

\textbf{Lemma 1}\hspace{0.5in}
\textit{In the equation (3.10), it follows that $\texttt{R}_{00}$ is divisible by $h_{00}$.}\\

($Proof$.)\hspace{0.5in}
Suppose that $\texttt{R}_{00}$ is not divisible by $h_{00}$ and since $(h_{ij})$ is positive definite, $(R_{00})^2$ 
is not divisible by $h_{00}$. 

Taking into account that ${P_{(5)}^i}_l$ and ${Q_{(9)}^i}_l$ are homogeneous polynomials of $y^i$ 
and that $(W_0)^2$ is  not divisible by $h_{00}$, it follows that the equation 
\begin{eqnarray*}
	h_{0l}y^i=h_{00}{\eta^i}_l,
\end{eqnarray*}
where ${\eta_l}^i(x)$ is a function of $(x^i)$ alone, holds good. Transvecting the above equation by $W^l$, we get
\begin{eqnarray*}
	W_0y^i=h_{00}{\eta_l}^i(x)W^l.
\end{eqnarray*}
Since $h_{00}$ is not divisible by $W_0$, the above equation is impossible. \hspace{0.5in}  Q.E.D.\\

Therefore, it follows that 
$\texttt{R}_{00}$ is divisible by $h_{00}$ and the following equation holds:
\begin{eqnarray*}
	\texttt{R}_{00}=c(x)h_{00},
\end{eqnarray*}
where $c(x)$ is a function of $(x^i)$ alone.
Derivating  the above equation by $y^i$ and $y^j$, we get
\begin{eqnarray}
	W_{i||j}+W_{j||i}=2c(x)h_{ij}.
\end{eqnarray}
Transvecting (3.11) by $W^iW^j$, we get $W_{i||j}W^iW^j=c(x)h_{ij}W^iW^j$ and using $h_{ij}W^iW^j=|W|^2=1$ and 
$W_{i||r}W^i=0$, we obtain $c(x)=0$.
Therefore, it follows that the equality
\begin{eqnarray}
	\texttt{R}_{ij}=0
\end{eqnarray}
holds good. Hence, we have that  $W$ is a Killing vector field.
Therefore, we can state\\

\textbf{Lemma 2}\hspace{0.5in}
\textit{If a Kropina space $(M, \alpha^2/\beta)$ is of constant curvature $K$, then\\
\hspace{0.5in}  (1) \hspace{0.1in} $W(x)$ is a Killing vector field,\\
and then\\
\hspace{0.5in}(2) \hspace{0.1in} Killing part $R^i_{(9)l}=0$.}\\

Using (3.12), the equation (3.10) reduces to
\begin{eqnarray}
	(h_{00})^2{P_{(5)}^i}_l+{Q_{(9)}^i}_l=0
\end{eqnarray}
and we have the following equalities:
\begin{eqnarray}
	W_{i||j}=\texttt{S}_{ij}, \hspace{0.1in}	
	\texttt{S}_j=W_{i||j}W^i=0,\hspace{0.1in}
	W_{0||j}=\texttt{S}_{0j},\hspace{0.1in}
	W_{i||0}=\texttt{S}_{i0},\hspace{0.1in}
	W_{0||0}=0.
\end{eqnarray}

\subsection{The vanishing part.}\label{}

In this subsection, we will show that the equality $Q_{(9)l}^i=0$ holds from the conclusion $\texttt{R}_{00}=0$
in the previous subsection.

 Using (3.12) and (3.14), the $A's$, $B's$, $C's$, $D's$, $E's$, $H's$ and $J's$  reduce  to
\begin{eqnarray*}
	&&A_{(1)}^i:=-2{W^i}_{||0},\hspace{0.1in}
	       B_{(1)l}^i=-2{W^i}_{||0||l},\hspace{0.1in}
                 B_{(21)l}^i=2{W^i}_{||0}W_{0||l}, \hspace{0.1in}
                     C_{(0)l}^i=-2{W^i}_{||l},\\
      &&  C_{(11)l}^i=2{W^i}_{||0}W_l,\hspace{0.1in}
               C_{(21)l}^i=-4{W^i}_{||0}h_{0l},\hspace{0.1in}
                      D_{(1)l}^i=-2{W^i}_{||l||0},\hspace{0.1in} 	\\
     &&D_{(21)l}^i =2{W^i}_{||0||0}W_l+2{W^i}_{||0}W_{l||0},\hspace{0.1in}
	          D_{(32)l}^i=-4{W^i}_{||0||0}h_{0l},\hspace{0.1in}
	                E_{(0)l}^i   =4{W^i}_{||r}{W^r}_{||l},\\
      &&E_{(11)l}^i=-4{W^i}_{||r}{W^r}_{||0}W_l,\hspace{0.1in}
              E_{(22)l}^i =8{W^i}_{||0}W_{0||l}
                            +8{W^i}_{||r} {W^r}_{||0}h_{0l},\\
      &&H_{(01)rl}^i =2{W^i}_{||r}W_l+2W_r{W^i}_{||l},\hspace{0.1in}
              H_{(11)rl}^i =-4W_r{W^i}_{||0}W_l,\hspace{0.1in}\\
      &&H_{(12)rl}^i =-4(h_{0r}{W^i}_{||l}+{W^i}_{||r}h_{0l}+{W^i}_{||0}h_{rl}),\hspace{0.1in}
          H_{(21)rl}^i=4(W_r{W^i}_{||0}h_{0l}+h_{0r}{W^i}_{||0}W_l),\\
       &&J_{(11)l}^i =-4{W^i}_{||r}{W^r}_{||0}W_l,\hspace{0.1in}
                      J_{(22)l}^i=8{W^i}_{||r}{W^r}_{||0}h_{0l}+8 {W^i}_{||0}W_{l||0},
 \end{eqnarray*}
and the others  are zero.
Using the above equalities, we get
\begin{eqnarray*}
      &&2J_{(31)l}^i -4D_{(31)l}^i-E_{(31)l}^i=0,\hspace{0.2in}
                     2J_{(41)l}^i -E_{(41)l}^i-4{D^i}_{(41)l}=0\\
      &&8B_{(3)}^i-4D_{(33)l}^i+2J_{(33)l}^i-E_{(33)l}^i=0,\hspace{0.2in}
                    2J_{(42)l}^i+8B_{(4)}^i-4 D_{(42)l}^i-E_{(42)l}^i=0,\\
      &&2J_{(51)l}^i-E_{(51)l}^i=0,\hspace{0.2in}
                    2J_{(43)l}^i-E_{(43)l}^i=0,\\
      &&2J_{(52)l}^i-4D_{(5)l}^i-E_{(52)l}^i=0,\hspace{0.2in}
	              2J_{(61)l}^i-4D_{(6)l}^i-E_{(61)l}^i=0.	
\end{eqnarray*}
Therefore, from Lemma 2, it follows \\

\textbf{Lemma 3}\hspace{0.5in}
\textit{If a Kropina space $(M, \alpha^2/\beta)$ is of constant curvature $K$, then\\
\hspace{0.5in}  (1) \hspace{0.1in}vanishing part $Q^i_{(9)l}=0$,\\
and then\\
\hspace{0.5in}(2) \hspace{0.1in} curvature  part $P^i_{(5)l}=0$.}\\

\subsection{The curvature part.}\label{}

\hspace{0.2in}
In this subsection, we will show that Lemma 3 implies that $(M, h)$ is a Riemannian space of constant curvature $K$.

 Using the resuls given at the beginning of the previous subsection, we have
\begin{eqnarray*}
	-E_{(0)l}^i-4K{\delta^i}_l&=&-4{W^i}_{||r}{W^r}_{||l}-4K{\delta^i}_l,\\
	2J_{(11)l}^i -E_{(11)l}^i-4KW_ly^i  &=&-4{W^i}_{||r}{W^r}_{||0}W_l-4KW_ly^i, \\
	8B_{(1)l}^i-4D_{(1)l}^i+2J_{(12)l}^i-E_{(12)l}^i&=&-16{W^i}_{||0||l}+8{W^i}_{||l||0},\\
     -4D_{(21)l}^i+8B_{(21)l}^i+2J_{(22)l}^i-E_{(22)l}^i+8Kh_{0l}y^i
	              &=&	-8{W^i}_{||0||0}W_l  +8{W^i}_{||r}{W^r}_{||0}h_{0l}   +8Kh_{0l}y^i,\\	
	16 {{^hR_0}^i}_{0l}+8B_{(22)l}^i-4D_{(22)l}^i+2J_{(23)l}^i-E_{(23)l}^i
	          &=&16 {{^hR_0}^i}_{0l},\\
	2J_{(32)l}^i-4D_{(32)l}^i -E_{(32)l}^i&=&16{W^i}_{||0||0}h_{0l}.
\end{eqnarray*}
Therefore, from (2) of Lemma 3 and the above  equalities we have
\begin{eqnarray}
	&&-\frac{1}{4}{P_{(5)}^i}_l=(h_{00})^2W_0({W^i}_{||r}{W^r}_{||l}+K{\delta^i}_l )
                         +(h_{00})^2({W^i}_{||r}{W^r}_{||0}W_l+KW_ly^i )\\
                  &&   +2h_{00}(W_0)^2(2{W^i}_{||0||l}-{W^i}_{||l||0})
                   +2h_{00}W_0({W^i}_{||0||0}W_l  -{W^i}_{||r}{W^r}_{||0}h_{0l}   -Kh_{0l}y^i)\nonumber\\
                  &&   -4(W_0)^3  \cdot {{^hR_0}^i}_{0l}
                       -4(W_0)^2{W^i}_{||0||0}h_{0l} =0,\nonumber
\end{eqnarray}

First, we consider the term $(h_{00})^2({W^i}_{||r}{W^r}_{||0}W_l+KW_ly^i )$ which does not contain $W_0$.
Taking into account  that $(h_{00})^2$ is not divisible by $W_0$, the following equality must hold good
\begin{eqnarray}
	 {W^i}_{||r}{W^r}_{||0}W_l+KW_ly^i=W_0{c_l}^i(x),
\end{eqnarray}
where ${c_l}^i(x)$ are functions of $(x^i)$ alone, must  hold.
Transvecting (3.16) by $y^l$ and dividing it by $W_0$, we get
\begin{eqnarray*}
	 {W^i}_{||r}{W^r}_{||0}+Ky^i={c_l}^i(x)y^l.
\end{eqnarray*}
Derivating it by $y^l$, we have
\begin{eqnarray*}
	 {W^i}_{||r}{W^r}_{||l}+K{\delta^i}_l={c_l}^i(x).
\end{eqnarray*}
Substituting the above equality in (3.16), we get
\begin{eqnarray*}
	 {W^i}_{||r}{W^r}_{||0}W_l+KW_ly^i=W_0({W^i}_{||r}{W^r}_{||l}+K{\delta^i}_l).
\end{eqnarray*}
Transvecting the above equality by $W^l$, we get
\begin{eqnarray*}
	 {W^i}_{||r}{W^r}_{||0}+Ky^i=KW_0W^i,
\end{eqnarray*}
where we have used ${W^i}_{||r}{W^r}_{||l}W^l=0$.
Derivating the above equality by $y^l$, we get
\begin{eqnarray}
	 {W^i}_{||r}{W^r}_{||l}=KW_lW^i-K{\delta^i}_l.
\end{eqnarray}
Substituting (3.17) in (3.15) and dividing it by $2W_0$, it follows
\begin{eqnarray}
	&&K(h_{00})^2W_lW^i   +h_{00}W_0(2{W^i}_{||0||l}-{W^i}_{||l||0}) \\
                  &&\hspace{0.6in}  +h_{00}({W^i}_{||0||0}W_l  -KW_0W^ih_{0l})
                                    -2(W_0)^2  \cdot {{^hR_0}^i}_{0l}
                                    -2W_0{W^i}_{||0||0}h_{0l} =0.\nonumber
\end{eqnarray}
Transvecting the above equality by $h_{0i}$, we get
\begin{eqnarray}
       {W_0}_{||l||0}=K(h_{00}W_l  -  W_0h_{0l}).                                                               
\end{eqnarray}
Using $W_{i||j}+W_{j||i}=0$ and (3.19), we have
\begin{eqnarray}
       {W_l}_{||0||0}=-K(h_{00}W_l  -  W_0h_{0l}).                                                               
\end{eqnarray}
Derivating (3.19) by $y^i$, we have
\begin{eqnarray*}
       {W_i}_{||l||0}+{W_0}_{||l||i} =K(2h_{0i}W_l  -  W_ih_{0l}-  W_0h_{il}).                                                               
\end{eqnarray*}
From the above equality, we have
\begin{eqnarray}
	{W_i}_{||0||l}  =-{W_i}_{||l||0}-K(2h_{0l}W_i  -  W_lh_{0i}-  W_0h_{li}).
\end{eqnarray}
Using (3.20) and (3.21), the equality (3.18) reduces to
\begin{eqnarray*}
	   3h_{00}(KW_ly^i -Kh_{0l}W^i  -{W^i}_{||l||0} )                
              -2W_0( {{^hR_0}^i}_{0l}   + K y^i h_{0l}- K h_{00}{\delta^i}_l)=0.
\end{eqnarray*}
Since $h_{00}$ is not divisible by $W_0$, it follows that the equality 
\begin{eqnarray}
	{{^hR_0}^i}_{0l}+K h_{0l}y^i-Kh_{00}{\delta^i}_l=h_{00}{d^i}_l(x),
\end{eqnarray}
where ${d^i}_l(x)$ are functions of $(x^i)$ alone, must hold.  Transvecting the above equality by $y^l$, we get
${d^i}_l(x)y^l=0$.
Derivating the above equality by $y^l$, we get ${d^i}_l(x)=0$. Substituting it in  (3.22), we get
\begin{eqnarray}
	{{^hR_0}^i}_{0l}+K h_{0l}y^i-Kh_{00}{\delta^i}_l=0.
\end{eqnarray}
We can rewrite the above equality as
\begin{eqnarray*}
	{{^hR_0}^i}_{0l}&=&K(h_{00}{\delta^i}_l- h_{0l}y^i).
\end{eqnarray*}
Derivating (3.23) by $y^j$ and $y^k$, we get
\begin{eqnarray}
	{{^hR_j}^i}_{kl}+{{^hR_k}^i}_{jl}=K(2h_{jk}{\delta^i}_l-h_{jl}{\delta^i}_k-h_{kl}{\delta^i}_j),
\end{eqnarray}
and interchanging $j$ and $l$, we obtain
\begin{eqnarray}
	{{^hR_l}^i}_{kj}+{{^hR_k}^i}_{lj}=K(2h_{lk}{\delta^i}_j-h_{lj}{\delta^i}_k-h_{kj}{\delta^i}_l).
\end{eqnarray}
Subtracting (3.25) from (3.24), we get
\begin{eqnarray*}
	{{^hR_j}^i}_{kl}+2{{^hR_k}^i}_{jl}-{{^hR_l}^i}_{kj}=3K(h_{jk}{\delta^i}_l-h_{kl}{\delta^i}_j).
\end{eqnarray*}
Since the left-hand side of the above equality can be changed as follows:
\begin{eqnarray*}
	{{^hR_j}^i}_{kl}+2{{^hR_k}^i}_{jl}-{{^hR_l}^i}_{kj}
	=2{{^hR_k}^i}_{jl}-{{^hR_j}^i}_{lk}-{{^hR_l}^i}_{kj}
	=2{{^hR_k}^i}_{jl}+{{^hR_k}^i}_{jl}
	=3{{^hR_k}^i}_{jl},
\end{eqnarray*}
we obtain
\begin{eqnarray*}
	{{^hR_k}^i}_{jl}=K(h_{jk}{\delta^i}_l-h_{kl}{\delta^i}_j).
\end{eqnarray*}
This means that the Riemannian space $(M, h)$ is of constant curvature $K$.\\

Therefore, we obtain \\

\textbf{Theorem 2}\hspace{0.5in}
 \textit{Let $M$ be an $n(\ge 2)$-dimensional Riemannian manifold. Put $\alpha=\sqrt{a_{ij}(x)y^iy^j}$ and 
$\beta=b_i(x)y^i$.
Let $(M, \alpha^2/\beta)$ be a Kropina space and define a new Riemannian metric $h=\sqrt{h_{ij}(x)y^iy^j}$ and
 a vector field $W$ with $|W|=1$ by (1.2) and (1.3).}

\textit{If the Kropina space $(M, \alpha^2/\beta)$ is of constant curvature $K$, then the vector field $W$ is a
 Killing one and the Riemannian space $(M, h)$ is of constant curvature $K$  .}

\subsection{The converse of Theorem 2.}\label{}

\hspace{0.2in}
Let $(M, \alpha^2/\beta)$ be a Kropina space and define a new Riemannian metric $h=\sqrt{h_{ij}(x)y^iy^j}$ and
 a vector field $W$ with $|W|=1$ by (1.2) and (1.3). 
Suppose that the vector field $W$ is a Killing one and that the Riemannian space $(M, h)$ is of constant curvature $K$.
To prove that  the Kropina space $(M, \alpha^2/\beta)$ is of constant curvature $K$, we  have only to show that 
the equality (3.9) holds.

Since the vector field $W$ is a Killing one, we have $\texttt{R}_{00}=0$. 
Taking into account (2) of Lemma 2 and (1) of Lemma 3, the Killing part $R$ and the vanishing part $Q$ vanishes 
respectively.

Therefore, we have only  to show that the curvature part $P^i_{(5)l}$ vanishes. At the rest of this subsection, we will prove it.
The curvature part $P^i_{(5)l}$ is defined by (3.15).

First, we give the following Lemma 4:\\

\textbf{Lemma 4}\hspace{0.5in}
\textit{ For a Killing vector field $W=W^i(\partial/\partial x^i)$ of constant length $|W|=1$, the equality 
\begin{eqnarray}
	W_{i||j||k}={W_r}\hspace{0.05in}{^h{{R_k}^r}_{ij}}
\end{eqnarray}
holds good.}\\

($Proof$.)\hspace{0.5in}
From the Ricci's formula, it follows
\begin{eqnarray}
	W_{i||j||k}-W_{i||k||j}=-{W_r}\hspace{0.05in} {^h{{R_i}^r}_{jk}}.
\end{eqnarray}

On the other hand, since $W$ is a Killing vector field, we have
\begin{eqnarray*}
	W_{i||j}+W_{j||i}=0.
\end{eqnarray*}
Applying the $h$-covariant derivative ${_{||k}}$ to the above equality, we get
\begin{eqnarray*}
	W_{i||j||k}+W_{j||i||k}=0.
\end{eqnarray*}
Replacing $i$, $j$, $k$ by $j$, $k$, $i$ and $k$, $i$, $j$ respectively, we obtain the following equalities:
\begin{eqnarray*}
	W_{j||k||i}+W_{k||j||i}&=&0,\\
       W_{k||i||j}+W_{i||k||j}&=&0.	
\end{eqnarray*}
Subtracting the second equality from the first equality and adding the third equality to it, we get
\begin{eqnarray}
	W_{i||j||k}+W_{i||k||j}-{W_r}\hspace{0.05in} {^h{{R_j}^r}_{ik}}-{W_r}\hspace{0.05in} {^h{{R_k}^r}_{ij}}=0.
\end{eqnarray}

From (3.27) and (3.28), we get
\begin{eqnarray*}
	2W_{i||j||k}=-{W_r}\hspace{0.05in} {^h{{R_i}^r}_{jk}}
	                                 +{W_r}\hspace{0.05in} {^h{{R_j}^r}_{ik}}
	                                         +{W_r}\hspace{0.05in} {^h{{R_k}^r}_{ij}}.
\end{eqnarray*}
Using the formula
\begin{eqnarray*}
	{^h{{R_i}^r}_{jk}}+{^h{{R_j}^r}_{ki}}+{^h{{R_k}^r}_{ij}}=0,
\end{eqnarray*}
we have
\begin{eqnarray*}
	2W_{i||j||k}&=&{W_r}\hspace{0.05in} {^h{{R_j}^r}_{ki}}
	                       +{W_r}\hspace{0.05in}{^h{{R_k}^r}_{ij}}
	                                 +{W_r}\hspace{0.05in} {^h{{R_j}^r}_{ik}}
	                                         +{W_r}\hspace{0.05in} {^h{{R_k}^r}_{ij}}\\
	           &=&2 {W_r}\hspace{0.05in}{^h{{R_k}^r}_{ij}}
	                                 	                                   \end{eqnarray*}
that is, (3.26) holds good.        \hspace{1in} Q.E.D.\\

From the assumption that the Riemannian space $(M, h)$ is of constant curvature, we have
\begin{eqnarray}
	{^h{{R_k}^r}_{ji}}=K(h_{kj}{\delta^r}_i-h_{ki}{\delta^r}_j).
\end{eqnarray}

Using the above equality and (3.26), we get
\begin{eqnarray}
	{W^i}_{||j||k}  =K({\delta^i}_kW_j-h_{kj}W^i)
\end{eqnarray}
and from here and   ${y^i}_{||j}=0$ (See, Remark 1), it follows 
\begin{eqnarray}
	{W^i}_{||0||l}&=&K({\delta^i}_lW_0-h_{l0}W^i),\nonumber\\
	{W^i}_{||0||0}&=&K(y^iW_0-h_{00}W^i),\\
	{W^i}_{||l||0}&=&K(y^iW_l-h_{0l}W^i).\nonumber	
\end{eqnarray}

From (3.29), we have
\begin{eqnarray}
	{^h{{R_0}^i}_{0l}}=K(h_{00}{\delta^i}_l-h_{0l}y^i)
\end{eqnarray}
and applying  the $h$-covariant derivative ${_{||i}}$ to the equality $|W|^2=W_rW_sh^{rs}=1$, we get
\begin{eqnarray*}
	W_{r||i}W^r=-W_{i||r}W^r=0.
\end{eqnarray*}
Furthermore, applying the $h$-covariant derivative ${_{||l}}$ to the above equality, we obtain 
\begin{eqnarray*}
	W_{i||r}{W^r}_{||l}+W_{i||r||l}W^r=0.
\end{eqnarray*}
From the above equality and (3.30), we have
\begin{eqnarray}
	W_{i||r}{W^r}_{||l}&=&-W_{i||r||l}W^r\\
	                   &=&K(h_{lr}W_i-h_{li}W_r)W^r  \nonumber\\
                         &=&K(W_lW_i-h_{li}).\nonumber
\end{eqnarray}
Substituting the equalities (3.31)-(3.33) in (3.15), we can easily recognize the curvature part ${P_{(5)}^i}_l=0$. 
 Therefore, (3.9) holds good. Hence, from Proposition 2, we get\\

\textbf{Theorem 3}\hspace{0.5in}
\textit{Let $M$ be an $n(\ge 2)$-dimensional Riemannian space. Put $\alpha=\sqrt{a_{ij}(x)y^iy^j}$ and $\beta=b_i(x)y^i$.
Let $(M, \alpha^2/\beta)$ be a Kropina space and define a new Riemannian metric $h=\sqrt{h_{ij}(x)y^iy^j}$ and
 a vector field $W=W^i(\partial/\partial x^i)$ of constant length $|W|=1$ by (1.2) and (1.3).}

\textit{If the vector field $W=W^i(\partial/\partial x^i)$ is a Killing one and the Riemannian space $(M, h)$ is of 
constant curvature $K$, the Kropina space $(M, \alpha^2/\beta)$ is of constant curvature $K$}.\\

From Theorem 2 and Theorem 3, we have\\

\textbf{Theorem 4}\hspace{0.5in}
\textit{Let $(M, \alpha^2/\beta)$ be an $n(\ge 2)$-dimensional Kropina space, where $\alpha^2=a_{ij}(x)y^iy^j$, 
$\beta=b_i(x)y^i$ and the matrix $(a_{ij})$ is positive definite. 
For the Kropina space, we define a new Riemannian metric $h=\sqrt{h_{ij}(x)y^iy^j}$ and a vector field 
$W=W^i(\partial/\partial x^i)$ of constant length $|W|=1$ on $M$ by (1.2) and (1.3).}

\textit{Then, the Kropina space $(M, \alpha^2/\beta)$ is  of constant curvature $K$ if and only if the following conditions
hold:}
 
\textit{(1)$W_{i||j}+W_{j||i}=0$, that is, $W=W^i(\partial/\partial x^i)$ is a Killing vector field.}

\textit{(2)The Riemannian space $(M, h)$ is of constant curvature $K$.}\\

Let $(M, F=\alpha^2/\beta)$ be an $n(\ge 2)$-dimensional Kropina space. From Theorem 1, for this Kropina metric  $F=\alpha^2/\beta$, we can
define  a Riemannian metric $h=\sqrt{h_{ij}(x)y^iy^j}$ and a  vector field $W=W^i(\partial/\partial x^i)$ of constant
length 1 on $M$ 
 by (1.2) and (1.3). We suppose that the vector field $W$ is a Killing one. Then, we have ${\texttt{R}}_{00}=0$. 
From this assumption, we get the second equation of (3.14), that is, ${\texttt{S}}_0=0$.
Substituting  ${\texttt{R}}_{00}=0$, ${\texttt{S}}_0=0$  and $F=h_{00}/(2W_0)$ in (2.6), we obtain the equation $\Phi^i=-F{\texttt{S}^i}_0$.
Substituing this in (2.4), we get\\

\textbf{Theorem 5}\hspace{0.5in}
\textit{Let $(M, F=\alpha^2/\beta)$ be an $n(\ge 2)$-dimensional  Kropina space. We define  a Riemannian metric 
$h=\sqrt{h_{ij}(x)y^iy^j}$ and a vector field $W=W^i(\partial/\partial x^i)$ of constant length $|W|=1$ on $(M,h)$  
by (1.2) and (1.3).}

\textit{Suppose that the vector field $W$ is a Killing one, then  the coefficients $G^i$ of the geodesic spray  of the
 Kropina space $(M, \alpha^2/\beta)$ is written as follows:
\begin{eqnarray}
	2G^i&=&{{^h\gamma_0}^i}_0-2F{\texttt{S}^i}_0,                      
\end{eqnarray}
where ${{^h\gamma_j}^i}_k$ are Christoffel symbols of the Riemannian space $(M, h)$.}\\

\textbf{Remark 2}\hspace{0.5in}
\textit{The geodesic spray of the Randers space $(M, \alpha+\beta)$ is given in the subsection 2.3  ([5] p.5).
 Comparing to this, the geodesic spray of the Kropina space $(M, F=\alpha^2/\beta)$ is in a very simple form (3.33).} \\

\textit{Acknowledgments.}\hspace{0.5in}
The authors would like to express their sincere thanks to Professor Dr. S. V. Sabau and Professor Dr. H. Shimada
for the valuable advices, the efforts to check all the calculation  and the continuous encouragements.

\vspace {0.5in}
\begin{center}
	References
\end{center}
\begin{description}	 
 \item[[1]] S. B\'acs\'o, X. Cheng and Z. Shen : Curvature properties of $(\alpha, \beta)$-metrics, 
                          Advanced Studies in Pure Mathematics 48, 2007, Finsler Geometry, Sappro 2005 -
                           In Memory of Makoto Matsumoto,  73-110.
\item[[2]] M. Matsumoto : Foundations of Finsler geometry and special Finsler spaces, Kaiseisha Press, Saikawa, Otsu,
       Japan, 1986.
\item[[3]] M. Matsumoto: Finsler spaces of constant curvature with Kropina metric, Tensor N.S., 50(1991), 194-201.
\item[[4]] M. Matsumoto: Theory of Finsler spaces with $(\alpha, \beta)-$metric, Rep. Math. Phys., 31(1992), 43-83.
\item[[5]] C. Robles : Geodesic in Randers spaces of constant curvature, arXiv:math/0501358v2 [math.DG] 8 July 2005.
\item [[6]] S. Tachibana : Riemannian geometry (Japanese), Modern mathematical lecture 15, Asakura Shoten,
                 Tokyo, 1967.
 \item[[7]] R. Yoshikawa and K. Okubo: Kropina spaces of constant curvature, Tensor N.S., 68(2007), 190-203.
\end{description}

\vspace{0.1in}
\begin{center}
	  \hspace{0.7in}  Hachiman technical          \hspace{1in}      Faculty of Education \\   
	 \hspace{0.8in}   High School                   \hspace{1.4in}    Shiga  University      \\   
	 \hspace{0.2in} 5 Nishinosho-cho Hachiman     \hspace{0.8in}      2-5-1 Hiratsu Otsu            \\   
	 \hspace{0.6in} 523-0816 Japan                \hspace{1.3in}    520-0862 Japan                     \\
	 \hspace{0.7in} E-mail: ryozo@e-omi.ne.jp     \hspace{0.4in}    E-mail: okubo@edu.shiga-u.ac.jp
\end{center}	

\end{document}